\documentclass[11pt]{article}
\usepackage{amsmath,amssymb,amsthm,epsfig}

\newtheorem{1}{Proposition}
\newtheorem{2}[1]{Proposition}
\newtheorem{3}[1]{Lemma}
\newtheorem{4}[1]{Theorem}
\newtheorem{5}[1]{Corollary}

\begin{document}

\title{Functions for relative maximization}
\author{E. Garibaldi e A. O. Lopes}
\date{\today}
\maketitle

\begin{abstract}

We introduce functions for relative maximization in a general
context: the beta and alpha applications. After a systematic study
concerning regularities, we investigate how to approximate certain
values of these functions using periodic orbits. We establish yet
that the differential of an alpha application dictates the
asymptotic behavior of the optimal trajectories.

\end{abstract}

\vspace {.8cm}

1. \textsc {{\large Introduction}}

\vspace {.4cm}

Let $ (X, d) $ be a compact metric space. If $ T: X \to X $ is a
continuous function, consider $ \mathcal M_T $ the set of the
$T$-invariant Borel probability measures. Remind that $\mathcal
M_T $ is convex and weak* compact.

Given a continuous function $ A: X \to \mathbb{R} $, we denote
$$ \beta_A = \max_{\mu \in \mathcal M_T} \int A \; d\mu. $$
In  ergodic optimization on compact spaces, the characterization
of the invariant probability measures whose integral of $ A $
reaches the maximum value $ \beta_A $ is one of the main goals. We
call any of these probabilities an $A$-maximizing probability.
References, general definitions and problems consider in this
theory can be found, for instance, in Jenkinson's notes (see
\cite{Jenkinson4}).

In the present work, we will look at the problem of the
description of maximizing probabilities in a slight different
formulation. Given $A$ as above, we introduce also a continuous
application $ \varphi: X \to \mathbb{R}^n $ (which plays the role
of a \emph{constraint}) and we extend the concept of $ \beta_A $
to a real function defined in a convex subset of $ \mathbb R^n $,
the  \emph{rotation set}. We use here a terminology borrowed from
the Aubry-Mather theory \cite{CI}. This concave application will
be called the \emph{beta function} $\beta_{A,\varphi}$ (associated
to $ A $ and $ \varphi $) and its Fenchel transform, the
\emph{alpha function}.

A motivation for analyzing such kind of problem is furnished by
\cite{LT3}. In this paper, fixing the Lagrangian $ L(x,v) = \| v
\|_x^2 / 2 $ obtained from the Riemannian metric in a compact
constant negative curvature surface of genus $ 2 $, the authors
focus Mather measures for $ L $ subject to a certain homological
condition. Via the Bowen-Series transformation $ T: S^1 \to S^1 $,
it is shown that this situation can be translated into a relative
maximization question for the potential $ A = - \log T' $.

A simple example of the theory studied here is the following.
Consider $ X=\{1,2,3\}^\mathbb{N} $ and $ A(x_0,x_1,x_2,\ldots) =
A(x_0,x_1) $ depending only on the two first coordinates. Suppose
$ A(2,3) = A(3,2) $ and $ A(2,3) > A(i,j)$ for the other
possibilities. The maximizing shift-invariant probability for $ A
$ without constraint is the one supported on the periodic orbit $
(2,3,2,3, \ldots) $. If we require $ {\displaystyle \int \varphi
\; d\mu = 1} $, where $ \varphi $ is the indicator function of the
cylinder $ \bar{1} $, then the maximization is given by the fixed
point $(1,1,\ldots)$.

In the setting presented here, we could change maximization for
minimization and the analogous statements would be immediately
verify. In the first section of the present paper, all the
definitions will be carefully discussed. We will then derive some
basic properties associated to these concepts.

The next topic will be dedicated to the study of the behavior of
alpha and beta applications when the constraint is modified. We
will show, for instance, the Lipschitz character of the
correspondence associating constraints to alpha functions values.
Concerning the influence over a beta function, we will verify
typically continuity of the respective correspondence.

We will also investigate the possibility of approximating alpha
and beta functions values using probability measures supported on
periodic orbits. Nevertheless, we will deal with this problem
specifically for the symbolic dynamics case. Under the hypothesis
of joint recurrence (to be discussed later), we will show the
existence of  periodic orbits carrying out the task.

We will conclude presenting a theorem that points out an
interesting connection between the differential of an alpha
application and the asymptotic behavior of certain trajectories
(see \cite{Gomes} for a different setting). In the proof of this
result, the sub-action concept will appear. This notion has been
largely studied \cite{Branco, Bousch1, Bousch2, CLT, Jenkinson3,
LT1, Savchenko, Souza}.

This paper is part of the first author's PhD thesis
\cite{Garibaldi}. It can be seen as an analysis of the properties
of the beta function $\beta_{A,\varphi}$ -- a generalization of
the maximal constant $ \beta_A $ -- and of its Fenchel transform.
On the other hand, once the graph of a beta application is part of
the boundary of a rotation set contained in $\mathbb R^{n + 1} $,
this study also brings some information on such set. However,
evidently it does not make it in a so explicit way as, for
instance, Kwapisz \cite{Kwapisz1, Kwapisz2, Kwapisz3} for certain
rotation sets arising from two-torus maps homotopic to the
identity, Bousch \cite{Bousch1} and Jenkinson \cite{Jenkinson2}
when analyzing the set of the barycentres of invariant measures
for circle maps or Ziemian \cite{Ziemian} in symbolic dynamics.

Finally we would like to point out that results presented here
should be seen as the general abstract setting of the theory. In
the future, we will analyze similar problems with more stringent
hypothesis, where the concept of non-crossing trajectories will be
present. We believe it will be possible to obtain stronger results
in this case.

\vspace {.8cm}

2. \textsc {{\large First Definitions}}

\vspace {.4cm}

Let $ \varphi: X \to \mathbb R^n $ be a continuous application
with coordinate functions $\varphi_1, \ldots, \varphi_n $. We have
then an induced map $\varphi_*: \mathcal M_T \to \mathbb R^n $
given by ${\displaystyle \varphi_*(\mu) = \left(\int \varphi_1 \;
d\mu, \ldots, \int \varphi_n \; d\mu \right)} $. Clearly,
$\varphi_* $ is a continuous and affine function.

We call $ \varphi_*(\mu) $ the rotation vector of the measure $
\mu \in \mathcal M_T $. (When $ n = 1 $, we will prefer the
expression rotation number of the measure.) Note that the image $
\varphi_*(\mathcal M_T) \subset \mathbb R^n $ is a convex compact
set, inheriting it of $ \mathcal M_T $. We denominate $
\varphi_*(\mathcal M_T) $ a rotation set. For $ h \in
\varphi_*(\mathcal M_T) $, the fiber $\varphi_{*} ^ {-1}(h) $ is
called the rotation class of $ h $. Also $\varphi_{*} ^ {-1}(h)
\subset \mathcal M_T $ is a convex compact set.

\begin{1}
For the induced map $\varphi_*: \mathcal M_T \to
\varphi_*(\mathcal M_T) $, we verify:

\vspace{.2cm}

{\raggedright (i) if the fiber $\varphi_{*} ^ {-1}(h) $ is a
singleton set containing an ergodic measure, then $ h $ is an
extremal point of $\varphi_*(\mathcal M_T) $;}

\vspace{.2cm}

{\raggedright (ii) if $ h $ is an extremal point of
$\varphi_*(\mathcal M_T) $, then the extremal points of
$\varphi_{*} ^ {-1}(h) $ are ergodic measures.}
\end{1}

Actually, this proposition is just a general version of results
presented in Jenkinson's PhD thesis (see lemmas 3.2 and 3.3 of
\cite{Jenkinson1}).

For $ A \in C^0(X) $, we define the beta function $ \beta_{A,
\varphi}: \varphi_*(\mathcal M_T) \to \mathbb R $ by
$$ \beta_{A, \varphi}(h) = \sup \left \{\int A \; d\mu: \mu \in \varphi_{*} ^ {-1}(h) \right \}. $$
In this context, we call the function $ \varphi $ a constraint and
the function $ A $ a potential. Important objects will be the
probability measures belonging to the  rotation class of $ h $
that, on such set, maximize the integral of the potential $ A $.
In clearer terms, consider the set
$$ \text{\Large $\mathit m$}_{A, \varphi}(h) = \left \{\mu \in \varphi_{*} ^ {-1}(h): \int A \; d\mu = \beta_{A, \varphi}(h) \right \}. $$
If $ \mu \in \text{\Large $\mathit m$}_{A, \varphi}(h) $, we say
that $\mu $ is an $(A, h)$-maximizing probability.

Since the rotation class of $ h $ is a compact set, it is easy to
prove that $ \text{\Large $\mathit m$}_{A, \varphi}(h) $ is a
nonempty compact set. It follows that $ \beta_{A, \varphi}:
\varphi_*(\mathcal M_T) \to \mathbb R $ is a concave application.
Moreover, since the correspondence $ {\displaystyle \mu \mapsto
\int A \; d\mu} $ is continuous, we conclude that $ \beta_{A,
\varphi} $ is continuous in the whole rotation set.

These properties of a beta application legitimate the definiton of
a concave function $ \alpha_{A, \varphi}: \mathbb R^n \to \mathbb
R $ via Fenchel transform
$$ \alpha_{A, \varphi}(c) = \min_{h \in \varphi_*(\mathcal M_T)} [\langle c,h \rangle - \beta_{A, \varphi}(h)]. $$
Such application is called the alpha function (associated to $ A $
and $ \varphi $).

It is interesting to examine the behaviors of the beta and alpha
applications when the parameters that define them are changed. For
instance, we can question how a potential modification affects a
beta function. Given $ h \in \varphi_* (\mathcal M_T) $, in a
natural way we obtain a function $ \beta_{\cdot, \varphi}(h):
C^0(X) \to \mathbb R $ that, to each potential $ A $, simply
associates the value $ \beta_{A, \varphi} (h) $. It is not
difficult to verify that this application is Lipschitz, with $
\text{Lip}(\beta_{\cdot, \varphi}(h)) \le 1 $.

A first consequence of this fact is the Lipschitz regularity of an
alpha function, with $ \text{Lip}(\alpha_{A, \varphi}) \le \|
\varphi \| _0 $. Indeed, since
$$ \alpha_{A, \varphi}(c) = - \max_{h \in \varphi_*(\mathcal M_T)} \beta_{A - \langle c,h \rangle, \varphi} (h), $$
if we take $ h' \in \varphi_*(\mathcal M_T) $ such that $
\alpha_{A, \varphi}(c') = - \beta_{A - \langle c',h' \rangle,
\varphi}(h') $, we have
\begin{eqnarray*}
\alpha_{A, \varphi}(c) - \alpha_{A, \varphi}(c') & \le & \beta_{A - \langle c',h' \rangle, \varphi}(h') - \beta_{A - \langle c,h' \rangle, \varphi}(h') \\
& \le & | \langle c - c',h' \rangle | \; \le \; \|\varphi\|_0 \|c
- c'\|.
\end{eqnarray*}

A second immediate consequence is the following version of the
Fenchel inequality
\begin{eqnarray*}
\beta_{A, \varphi}(h) + \alpha_{B, \varphi}(c) & \le & \beta_{A, \varphi}(h) + \langle c,h \rangle -  \beta_{B, \varphi}(h) \\
& \le & \langle c,h \rangle + \| A - B \|_0.
\end{eqnarray*}
Using this inequality, we see that it is also Lipschitz the
application that makes to correspond $ \alpha_{A, \varphi}(c) $ to
each potential $ A $, namely, the function $ \alpha_{\cdot,
\varphi}(c): C^0(X) \to \mathbb R $, besides we have
$\text{Lip}(\alpha_{\cdot, \varphi}(c)) \le 1 $.

Some properties of the applications $ \beta_{\cdot, \varphi}(h),
\alpha_{\cdot, \varphi}(c): C^0(X) \to \mathbb R$ are summarized
in the proposition below. The simple proof will be omitted.

\begin{2}
If $ A, B \in C^0(X) $, $ a \in \mathbb R $ and $ t, t' \in [0, 1]
$ with $ t + t' = 1 $, then the functions $ \beta_{\cdot,
\varphi}(h), \alpha_{\cdot, \varphi}(c): C^0(X) \to \mathbb R$
verify

\vspace{.2cm}

{\raggedright (i) $ \beta_{a A, \varphi}(h) = | a |
\beta_{\text{sgn}(a)A, \varphi}(h) $;}

\vspace{.2cm}

{\raggedright (ii) $ \beta_{A + B \circ T - B + a, \varphi}(h) =
\beta_{A, \varphi}(h) + a $, $ \alpha_{A + B \circ T - B + a,
\varphi}(c) = \alpha_{A, \varphi}(c) - a $;}

\vspace{.2cm}

{\raggedright (iii) $ \beta_{A + B, \varphi}(h) \le \beta_{A,
\varphi}(h) + \beta_{B, \varphi}(h) $;}

\vspace{.2cm}

{\raggedright (iv) $ \beta_{tA + t'B, \varphi}(h) \le t\beta_{A,
\varphi}(h) + t'\beta_{B, \varphi}(h) $, $ \alpha_{tA + t'B,
\varphi}(c) \ge t\alpha_{A, \varphi}(c) + t'\alpha_{B, \varphi}(c)
$;}

\vspace{.2cm}

{\raggedright (v) $ A \le B $ implies $ \beta_{A, \varphi}(h) \le
\beta_{B, \varphi}(h) $, $ \alpha_{A, \varphi}(c) \ge \alpha_{B,
\varphi}(c) $.}
\end{2}

Note that complementing expressions for the items (i) and (iii)
would be
$$ \alpha_{a A, \varphi}(c) = | a | \alpha_{\text{sgn}(a)A, \varphi}(c / | a |)
\text{ (for $ a \ne 0 $) and } $$
$$ \alpha_{A + B, \varphi}(c + c') \ge \alpha_{A, \varphi}(c) + \alpha_{B, \varphi}(c'), $$
which are not properties of the application $ \alpha_{\cdot,
\varphi}(c) $.

Besides redefining a beta function, the modification of the
potential still redescribes the set of maximizing probabilities.
Though, a particularity typically prevails.

\begin{2}
Assume $ h \in \varphi_* (\mathcal M_T) $. There is a residual
subset $\mathcal G = \mathcal G (h) \subset C^0(X) $ such that,
for each potential $ A \in \mathcal G $, $ \text{\Large $\mathit
m$}_{A, \varphi}(h) $ contains an unique probability measure.
\end{2}

This result can be seen as a particular version of a more general
formulation obtained in proposition 10 of \cite{CLT}. The proof
there is also valid for any compact metric space $ X $.

\vspace {.8cm}

3. \textsc {{\large The role of the constraint}}

\vspace {.4cm}

Our objective now will be to discuss how the changing of the
constraint affects the beta and alpha functions.

Suppose $ A, B \in C^0(X) $ and $ \varphi, \psi \in C^0(X, \mathbb
R^n) $. Take yet $ a \in \mathbb R^* $, $ b \in \mathbb R^n $ and
$ t, t' \in [0, 1] $ with $ t + t' = 1 $. Then, we have

\vspace {.2cm}

(i) $ \beta_{A, a\varphi} (h) = \beta_{A, \varphi} (h/a) $, $
\alpha_{A, a\varphi} (c) = \alpha_{A, \varphi} (ac) $;

(ii) $ \beta_{A, \varphi + \psi \circ T - \psi + b} (h) =
\beta_{A, \varphi} (h - b) $, $ \alpha_{A, \varphi + \psi \circ T
- \psi + b} (c) = \alpha_{A, \varphi} (c) + \langle c,b \rangle $;

(iii) $ \alpha_{A + B, \varphi + \psi} (c) \ge \alpha_{A, \varphi}
(c) + \alpha_{B, \psi} (c) $;

(iv) $ \alpha_{tA + t'B, t\varphi + t'\psi} (c) \ge t \alpha_{A,
\varphi} (c) + t' \alpha_{B, \psi} (c) $;

(v) $ \text{\Large $\mathit m $} _{A, \varphi}(h) \cap
\text{\Large $\mathit m $}_{A, \psi}(h') \ne \emptyset \Rightarrow
t\beta_{A, \varphi}(h) + t'\beta_{A, \psi}(h') \le \beta_{A,
\varphi + \psi}(h + h') $.

\vspace {.2cm}

{\raggedright The proof of these items is left to the reader.}

In order to make interesting the investigation of the relationship
between constraints and the beta and alpha functions, notice that
there is an initial difficulty: the constraint also determines the
domain of a beta application. Therefore, we first need to
establish which effect the change of this parameter produces on
the rotation set.

For a complete metric space $ Y $, we will denote $ \mathcal K(Y)
$ the collection of its compact subsets. With the Hausdorff
metric, $ \mathcal K(Y) $ becomes a complete metric space. For the
proof, we need to note the linear operator $ *: C^0 (X, \mathbb
R^n) \to C^0 (\mathcal M_T, \mathbb R^n) $ is bounded, with norm
smaller or equal to 1.

\begin{2}
If the application $ \Gamma_T: C^0 (X, \mathbb R^n) \to \mathcal
K(\mathbb R^n) $ is given by $ \Gamma_T (\varphi) = \varphi_*
(\mathcal M_T) $, then $ \Gamma_T $ is Lipschitz, with $
\text{Lip}(\Gamma_T) \le 1 $.
\end{2}

\begin{proof}
Observe that, for any $ \varphi, \psi \in C^0 (X, \mathbb R^n) $
and $ \mu \in \mathcal M_T $, we have
\begin{eqnarray*}
d(\varphi_* (\mu), \psi_* (\mathcal M_T))
& = & \inf_{\nu \in \mathcal M_T} \| \varphi_*(\mu) - \psi_*(\nu) \| \\
& \le & \| \varphi_*(\mu) - \psi_*(\mu) \| \\
& \le & \| (\varphi - \psi)_* \|_0 \\
& \le & \| \varphi - \psi \|_0 .
\end{eqnarray*}
However, by the construction of the Hausdorf metric, this argument
is enough to establish the proposition.
\end{proof}

Suppose we associate, to each map $ T \in C^0(X, X) $, some
application $ \varphi_T \in C^0(X, \mathbb R^n) $. This happens,
for instance, in the works on rotation sets arising from
$n$-dimensional torus homeomorphisms homotopic to the identity or
when one wants to analyze the spectrum of the Lyapunov exponents
of a differential application. Motivated by the proposition above,
we could question which regularity the map $ T \mapsto
(\varphi_T)_*(\mathcal M_T) $ presents. The proposition 5 answers
this demand.

\begin{2}
Consider a subset $ \mathcal U \subset C^0(X, X) $ with the
induced topology. Let $ T \in \mathcal U \mapsto \varphi_T \in
C^0(X, \mathbb R^n) $ be a continuous map. Then the application $
\Gamma_{\mathcal U}: \mathcal U \to \mathcal K(\mathbb R^n) $
defined as $ \Gamma_{\mathcal U} (T) = (\varphi_T)_* (\mathcal
M_T) $ is upper semi-continuous.
\end{2}

\begin{proof}
If the upper semi-continuity of $\Gamma_{\mathcal U} $ was not
verified, this would mean the existence of a map $ T \in \mathcal
U $ and some $ \epsilon > 0 $ for which we could determine a
sequence $ \{T_j \} \subset \mathcal U $ convergent to $ T $ and a
sequence of Borel probability measures $ \{\mu_j \} $ satisfying
both $ \mu_j \in \mathcal M_{T_j} $ and $ d((\varphi_{T_j})
_*(\mu_j), \Gamma_{\mathcal U} (T)) \ge \epsilon $. However, we
would have a subsequence $ \{\mu_{j_k} \} $ convergent to some
Borel probability measure $ \mu $. Then, the possibility of, for
any function $ f \in C^0(X) $, passing to the limit in
$$ \int f \circ T_{j_k} \; d\mu_{j_k} = \int f \; d\mu_{j_k} $$
would bring a contradiction: $ \mu \in \mathcal M_T $ and $
d((\varphi_T)_*(\mu),(\varphi_T)_*(\mathcal M_T)) \ge \epsilon $.
\end{proof}

For rotation sets arising from $n$-dimensional torus continuous
maps homotopic to the identity, a result in the spirit of the
previous proposition was demonstrated by Misiurewicz and Ziemian
(see the theorem 2.10 of \cite{MZ}) and by Herman (see section 10
of the chapter 1 of \cite{Herman}).

We continue the search for the influences of the constraint on the
beta and alpha functions. The next lemma will be useful.

\begin{3}
Given a constraint $ \varphi \in C^0 (X, \mathbb R^n) $, take $ h
\in \varphi_*(\mathcal M_T) $. Consider a sequence of constraints
$ \{\varphi_j \} $ converging to $ \varphi $. It follows

\vspace{.2cm}

{\raggedright (i) $ {\displaystyle \lim_{j \to \infty} d(h,
(\varphi_j)_*(\text{\Large $\mathit m$}_{A, \varphi}(h))) = 0} $;}

\vspace{.2cm}

{\raggedright (ii) $ {\displaystyle \lim_{j \to \infty} d(h,
\varphi_*(\text{\Large $\mathit m$}_{A, \varphi_j}(h))) = 0} $
when $ h \in (\varphi_j)_* (\mathcal M_T) $;}

\vspace{.2cm}

{\raggedright (iii) $ {\displaystyle \lim_{j \to \infty} d(h,
\varphi_*(\text{\Large $\mathit m$}_{A, \varphi_j}(h_j))) = 0} $
when $ h_j \in (\varphi_j)_*(\varphi_*^{-1} (h)) $.}
\end{3}

\begin{proof}
By a similar reasoning to the one used in proposition 4, we obtain
$ d(h', \psi_*(\text{\Large$\mathit m $}_{A, \psi'} (h'))) \le \|
\psi - \psi' \|_0 $, which clearly gives the items (i) and (ii).

Besides, we have $ d(h, \varphi_*(\text{\Large $\mathit m $}_{A,
\varphi_j}(h_j))) \le \| h - h_j \| + d(h_j,
\varphi_*(\text{\Large $\mathit m $}_{A, \varphi_j}(h_j))) $.
Consequently, for the second parcel, just as in the previous
paragraph, we have $ d(h_j, \varphi_*(\text{\Large$\mathit m
$}_{A, \varphi_j}(h_j))) \le \| \varphi_j - \varphi \|_0 $. And
for the first term, choosing $ \mu \in \varphi_*^{-1}(h) \cap
(\varphi_j)_*^{-1}(h_j) $, we have $ \| h - h_j \| = \|
\varphi_*(\mu) - (\varphi_j)_*(\mu) \| \le \| \varphi - \varphi_j
\|_0 $, which concludes the proof of the item (iii).
\end{proof}

Using the lemma above, we can show a kind of \emph{holographic
continuity} for the beta application as function of the
constraint.

\begin{2}
About the behavior of the beta and alpha functions when the
constraint is modified, we have the following results.

\vspace {.2cm}

{\raggedright (I) Given a constraint $ \varphi \in C^0 (X, \mathbb
R^n) $, take $ h \in \varphi_*(\mathcal M_T) $. Let $\{\varphi_j
\} $ be a sequence of constraints convergent to $ \varphi $.
Assume that $ \{ h_j \} \subset \mathbb R^n $ is a sequence
satisfying $ h_j \in (\varphi_j)_*(\text{\Large $\mathit m $}_{A,
\varphi}(h)) $. Then $ \lim \beta_{A, \varphi_j} (h_j) = \beta_{A,
\varphi} (h) $.}

\vspace {.2cm}

{\raggedright (II) For $ c \in \mathbb R^n $, the map $ \varphi
\mapsto \alpha_{A, \varphi}(c) $ is Lipschitz, with $
\text{Lip}(\alpha_{A, \cdot}(c)) \le \| c \| $.}
\end{2}

\begin{proof}
(I) Initially, note that, from the choice of $ h_j $, it happens $
\beta_{A, \varphi} (h) \le \beta_{A, \varphi_j} (h_j) $. Define a
sequence $ \{\eta_j \} \subset \mathbb R^n $ such that, for each
integer $ j $, the vector $ \eta_j \in \varphi_*(\text{\Large
$\mathit m $}_{A, \varphi_j}(h_j)) $ satisfies $\| h - \eta_j \| =
d(h,\varphi_*(\text{\Large $\mathit m $}_{A, \varphi_j}(h_j))) $.
Therefore, we obtain $ \beta_{A, \varphi} (h) \le \beta_{A,
\varphi_j} (h_j) \le \beta_{A, \varphi} (\eta_j) $. Besides, by
the item (iii) of the lemma above, we use $ \lim \eta_j = h $ to
legitimate $ \lim \beta_{A, \varphi}(\eta_j) = \beta_{A, \varphi}
(h) $.

\vspace {.1cm}

(II) Given $ \epsilon > 0 $, consider $ h \in \varphi_*(\mathcal
M_T) $ accomplishing $ \langle c,h \rangle - \beta_{A, \varphi}(h)
< \alpha_{A, \varphi}(c) + \epsilon/2 $. Afterwards, take a
probability measure $ \mu \in \varphi_*^{-1}(h) $ satisfying $
{\displaystyle \int A \; d\mu} > \beta_{A, \varphi}(h) -
\epsilon/2 $. Therefore, $ \langle c, \varphi_*(\mu) \rangle -
{\displaystyle \int A \; d\mu}
 < \alpha_{A, \varphi}(c) + \epsilon $. Besides, if $ \psi \in C^0(X, \mathbb R^n) $
is a constraint, the Fenchel inequality gives $ \alpha_{A,
\psi}(c) + {\displaystyle \int A \; d\mu} \le \alpha_{A, \psi}(c)
+ \beta_{A, \psi}(\psi_*(\mu)) \le \langle c, \psi_*(\mu) \rangle
$. Thus, we verify $ \alpha_{A, \psi}(c) - \alpha_{A, \varphi}(c)
< \langle c, (\psi - \varphi)_*(\mu) \rangle + \epsilon \le \| c
\| \| \psi - \varphi\|_0 + \epsilon $. And the result follows from
the symmetrical role carried out by $ \varphi $ and $ \psi $ and
from the arbitrariness of $ \epsilon $. \end{proof}

The presence of the sequence $ \{h_j \} $ in the proposition 7.I
is a little bit disappointing. Because of the item (i) of  lemma 6
we can choose it converging to $ h $, but we should ask: when $
\lim \beta_{A, \varphi_j} (h) = \beta_{A,\varphi} (h) $? The first
aspect to be noted is the requirement $ h \in (\varphi_j)_*
(\mathcal M_T) $. However, if $ h \in \text{int}(\varphi_*
(\mathcal M_T)) $, the proposition 4 assures that, for a
constraint $ \psi $ sufficiently close to $ \varphi $, we have $ h
\in \text{int}(\psi_* (\mathcal M_T)) $.
\begin{2}
Let $ \varphi \in C^0 (X, \mathbb R^n) $ be a constraint. Take $ h
\in \text{int} (\varphi_*(\mathcal M_T)) $. If $\{ \varphi_j \} $
is a sequence convergent to $ \varphi $, then $ \lim \beta_{A,
\varphi_j} (h) = \beta_{A, \varphi} (h) $.
\end{2}

\begin{proof}
Without loss of generality, we can suppose that $ h \in
\text{int}((\varphi_j)_* (\mathcal M_T)) $. However, we will need
a stronger version of this supposition. Fortunately, the
proposition 4 also allows to assume that $ D_\epsilon [h] \subset
\text{int}((\varphi_j)_*(\mathcal M_T)) $, where $ D_\epsilon [h]
$ is a closed ball of center $ h $ and radius $ \epsilon > 0 $
contained in $ \text{int}(\varphi_* (\mathcal M_T)) $.

Define a sequence of probability measures $\{\mu_j \} \subset
\text{\Large $\mathit m $}_{A, \varphi}(h) $ such that, for each
integer $ j $, we have $\| h - (\varphi_j)_*(\mu_j) \| = d(h,
(\varphi_j)_*(\text{\Large$\mathit m $} _{A, \varphi}(h))) $.
Putting $ h_j = (\varphi_j)_*(\mu_j) $, we set
$$\epsilon_j = \frac {\| h - h_j \|} {\| h - h_j \| + {\displaystyle\frac{\epsilon}{3}}}. $$
Write, then, $ h'_j = h_j + \epsilon_j^{-1}(h - h_j) $. Observe
that, in reason of the item (i) of the lemma 6, for an integer $ j
$ sufficiently large, it happens $ \| h - h_j \| \le \epsilon/3 $.
Hence, for such indexes, we verify $ h'_j \in D_\epsilon [h]
\subset \text{int}((\varphi_j)_* (\mathcal M_T)) $, that is, we
obtain $ h'_j = (\varphi_j)_* (\mu'_j) $ for some $T$-invariant
probability measure $\mu'_j $.

Put, for $ j $ sufficiently large, $ \mu '' _j = \epsilon_j \mu'_j
+ (1 - \epsilon_j) \mu_j $. Note that $ (\varphi_j)_*(\mu '' _j) =
\epsilon_j h'_j + (1 - \epsilon_j) h_j = h $. Therefore, if the
vector $ \eta_j \in \varphi_*(\text{\Large $\mathit m $}_{A,
\varphi_j}(h)) $ accomplishes $\| h - \eta_j \| = d(h,
\varphi_*(\text{\Large $\mathit m $} _{A, \varphi_j}(h))) $, we
have
$$ \epsilon_j \int A \; d\mu'_j + (1 - \epsilon_j) \beta_{A, \varphi} (h) = \int A \; d\mu '' _j \le \beta_{A, \varphi_j} (h) \le \beta_{A,\varphi} (\eta_j). $$
And the result follows directly of the items (i) and (ii) of the
lemma 6.
\end{proof}

The proposition above admits a more direct demonstration, but
maybe less instructive. Actually, it would simply be enough to
apply the conclusion of the proposition 4 to the functions $ \Phi
= (\varphi, A) $ and $ \Phi_j = (\varphi_j, A) $. This argument
will be explored ahead in the text\footnote{See, for instance, the
proof of the proposition 14.}.

\vspace {.8cm}

4. \textsc {{\large Approximation by Periodic Orbits}}

\vspace {.4cm}

Although, in this section, we will restrict the class of dynamical
systems to be examined, limiting us to study the approximation
problem for periodic orbits in the context of the symbolic
dynamics, we draw a general itinerary in certain aspects. This
itinerary describes how, when the purpose is to estimate certain
values of a beta application or of an alpha function, we can find
probability measures supported on periodic orbits accomplishing
such task.

Some comments on definitions and notations will be useful. Note
that, in any probability space $(Y, \mathcal B, \nu) $, given an
integrable application $ f: Y \to \mathbb R^n $, we still have the
idea of rotation vector of the measure $ \nu $. The integrability,
in fact, plays the main role when we write $ {\displaystyle
f_*(\nu) = \left(\int f_1 \; d\nu, \ldots, \int f_n \; d\nu
\right)} $. Given an ergodic function $ F: Y \to Y $, we set $
\text{\Large $\mathit b $} (f) $ to indicate the set of the
elements of $ Y $ that, for the application $ f \in L^1(Y,
\mathcal B, \nu) $, satisfy the Birkhoff's ergodic theorem. For
the characteristic function of a measurable set $ D $, however, we
will prefer to denote it $ \text{\Large $\mathit b $} (D) $.
Besides, just looking at the mensurability of $ F $, we put $
{\displaystyle S_k f = \sum_{j = 0}^{k - 1} f \circ F^j} $ for $ k
> 0 $ and $ S_0 f = 0 $.

We consider $ \Xi (D) $ the set of the elements $ z $ of $ D $
such that, for any $ \epsilon > 0 $, it exists a positive integer
$ L $ accomplishing both $ F^L (z) \in D $ and $ \left \| S_L f(z)
- L f_*(\nu) \right \| < \epsilon $. Thus, we say that the
integrable function $ f $ is joint recurrent (in relation to the
probability measure $ \nu $) if, for each $ D \in \mathcal B $, it
happens $ \nu(\Xi (D)) = \nu(D) $. (When $ n = 1 $, we will simply
say that $ f $ is recurrent.) If we want to identify functions
verifying such property, the following proposition describes a
sufficient condition.

\begin{2}
Let $ (Y, \mathcal B, \nu) $ be a probability space. Consider an
ergodic transformation $ F: Y \to Y $ and an integrable function $
f: Y \to \mathbb R^n $ satisfying
$$ \lim_{k \to \infty} \frac{1}{k^{1/n}} \left \| S_k f(y) - kf_*(\nu) \right \| = 0 $$
for $\nu$-almost every point $ y \in Y $. Then $ f $ is joint
recurrent.
\end{2}

Note that, if $ n = 1 $, for every integrable application, we have
immediately the required limit by the Birkhoff's ergodic theorem.
In simple terms, the proposition 9 shows that any integrable
function $ f: Y \to \mathbb R $ is recurrent. This result when $ n
= 1 $ was used by Ma\~n\'e in one of his works on minimizing
measures of Lagrangian systems (consult the lemma 2.2 of
\cite{Mane}). Nevertheless, two decades before, a theorem
containing the one-dimensional version of proposition 9 had been
obtained by Atkinson in \cite{Atkinson}. The proof that we will
present for the general case $ n \ge 1 $ is a generalization of a
proof for the particular situation when $ n = 1 $, more
specifically, of the proof given for the lemma 3-6.4 of \cite{CI}.

\begin{proof}
Without loss of generality, we can take $ f_*(\nu) = 0 $. Suppose
$ D \in \mathcal B $ with $ \nu(D) > 0 $. Assuming $ \epsilon > 0
$, let $ \Xi_\epsilon (D) $ denote the set of points $ z \in D $
for which there is a positive integer $ L $ such that $ F^L(z) \in
D $ and $ \| S_L f(z) \| < \epsilon $. Since $ \Xi (D) = \bigcap
\Xi_{1/j} (D) $, it is enough to show that $ \nu (\Xi_\epsilon
(D)) = \nu (D) $.

Take $ y \in D \cap \text{\Large $\mathit b $} (f) \cap
\text{\Large $\mathit b $} (D) \cap \text{\Large $\mathit b $}
(\Xi_\epsilon (D)) $ such that $ \lim k^{-1/n} \left \| S_k f(y)
\right \| = 0 $. Consequently, let $ L_1 < L_2 < \ldots L_k <
\ldots $ be the positive integers such that $ F^{L_k}(y) \in D $.
Defining $ a_k = S_{L_k} f(y) $, consider yet
$$ R = \{k: \forall \; \; m > k, \; \; \| a_m - a_k \| \ge \epsilon \} \; \text{ and } \; R_k = R \cap \{1,\ldots k \}. $$

Notice then, for each $ l \in \{1,\ldots k \} - R_k $, there
exists $ m > l $ such that $ \| S_{L_m - L_l} f(F^{L_l}(y)) \| =
\| a_m - a_l \| < \epsilon $. In other words, $ l \in \{1,\ldots k
\} - R_k $ implicates $ F^{L_l}(y) \in \Xi_\epsilon (D) $.
Therefore, we verify
\begin{eqnarray*}
1 + \#R_k
& \ge & 1 + \#\{1 \le l < k: F^{L_l}(y) \notin \Xi_\epsilon (D) \} \\
& \ge & \#\{0 \le j < L_k: F^j(y) \in D - \Xi_\epsilon (D) \}.
\end{eqnarray*}
Hence, since
$$ \nu(D - \Xi_\epsilon (D)) = \lim_{k \to \infty} \frac{1}{L_k} \sum_{j = 0}^{L_k - 1} \chi_{D - \Xi_\epsilon (D)} (F^j(y)), $$
the proposition will be demonstrated when we obtain a subsequence
of $ \{(\#R_k)/L_k \} $ converging to zero.

If $ R $ is a finite set, there is nothing to argue. Suppose,
otherwise, $ R $ is an infinite set. Then, by construction, $
\{a_k: k \in R \} $ is unbounded. In such case, choose an infinite
sequence $ S \subset R $ such that, for every $ k \in S $,
$$\| a_k \| = \max_{l \in R_k} \| a_l \|. $$
When denoting $ D_\rho(\gamma) $ the open ball of center $ \gamma
\in \mathbb R^n $ and radius $ \rho > 0 $, we observe that, given
$ k \in S $, $ D_{\epsilon/2}(a_l) \subset D_{\| a_k \| +
\epsilon/2}(0) $ for each $ l \in R_k $. Besides, for the
definition of $ R $, these balls $ D_{\epsilon/2}(a_l) $, $ l \in
R_k $, are all disjoint. Consequently,
$$ \# R_k \le \frac{{\displaystyle \left ( \| a_k \| + \frac{\epsilon}{2} \right )^n}}{{\displaystyle \left ( \frac{\epsilon}{2} \right )^n}}
= \sum_{j = 0}^n \binom{n}{j} \| S_{L_k} f(y) \|^j \left (
\frac{2}{\epsilon} \right )^j . $$ Reminding that $ \lim k^{-1/n}
\left \| S_k f(y) \right \| = 0 $, to verify
$$ \lim_{k \in S} \frac {\#R_k}{L_k} = 0 $$
is a current task.
\end{proof}

Let $ f $  be a joint recurrent function in relation to a
probability measure $ \nu $. If $ D $ is a measurable set of
positive measure, write $ \Xi^j (D) = \Xi (\Xi^{j - 1} (D)) $.
Then, observe that $ \nu (\bigcap \Xi^j (D)) = \nu (D) > 0 $. In
particular, if we have $ E \in \mathcal B $ with $ \nu(E) = 1 $,
then $ \bigcap \Xi^j (D) \cap E \neq \emptyset $. This simple fact
will play a crucial role in the proof of the next result. We will
need, however, more structure to obtain the statement of next
theorem. Thus, our study will be driven towards the symbolic
dynamics setting.

Let us begin, nevertheless, reminding concepts no restricted to
this dynamics. Given a periodic point $ x \in X $ of period $ M $,
naturally we have a $T$-invariant probability measure defined by
$$ \mu = \frac{1} {\#\text {orb}(x)} \sum_{y \; \in \text {orb}(x)} {\delta_y} = \frac{1}{M} \sum_{j = 0}^{M - 1}{\delta_{T^j(x)}}. $$
A way to refer to a such measure $ \mu $ will be calling it a
periodic probability measure. When taking $ x, y \in X $ and any
positive interger $ k $, other item to be remembered is the
synthesis between the metric and the dynamics indicated by
$$ d_k (x, y) = \max_{0 \le j < k} d(T^j(x), T^j(y)). $$
A special collection of potentials will be the focus of  our work:
the Walters potentials. A function $ f \in C^0(X) $ is a Walters
function if it admits a Walters module, that is, if there exists a
function $ H: \mathbb R^+ \to \mathbb R^+ \cup \{+ \infty \} $
increasing, null and continuous in zero, such that
$$\forall \; s \in \mathbb R^+, \; \forall \; k > 0, \; \forall \; x, y \in X, \; d_k (x, y) \le s \Rightarrow | S_k f(x) - S_k f(y) | \le H(s). $$
In the construction of sub-actions and in the search for
maximizing measures, this regularity condition was introduced by
Bousch in \cite{Bousch2}. For hyperbolic dynamical systems, the
set of the Walters functions includes (see the
definition-proposition 2 of \cite{Bousch2}) all the functions of
summable variation, in particular the H\"older functions are then
examples of Walters functions.

Finally, let $ \sigma: \Sigma \to \Sigma $ be a subshift of finite
type. Given a constant $ \lambda \in (0, 1) $, we consider in $
\Sigma $ the metric $ d (\mathbf x, \mathbf y) = \lambda^k $,
where $ \mathbf x, \mathbf y \in \Sigma $, $\mathbf x = (x_0, x_1,
\ldots) $, $ \mathbf y = (y_0, y_1, \ldots) $ and $ k = \min \{j:
x_j \ne y_j \} $. We will say that a continuous function $ g:
\Sigma \to \mathbb R^n $ is locally constant if there exists an
integer $ j \ge 0 $ such that $ g(\mathbf x) = g(\mathbf y) $
whenever $ x_0 = y_0, \ldots, x_j = y_j $. We could also say that
this application depends on $ j + 1 $ coordinates.

\begin{4}
Suppose $ \varphi \in C^0(\Sigma, \mathbb Q^n) $ is a locally
constant constraint and $ A $ is a Walters potential. Let $\varphi
$ be a joint recurrent application in relation to an ergodic
probability measure $ \nu \in \varphi_*^{-1} (r) $, where $ r \in
\varphi_*(\mathcal M_\sigma) \cap \mathbb Q^n $. Then, for each
$\epsilon > 0 $, there exists a periodic probability measure $ \mu
\in \varphi_*^{-1} (r) $ such that $ {\displaystyle \left | \int A
\; d\nu - \int A \; d\mu \right | < \epsilon} $.
\end{4}

\begin{proof}
Take $\mathbf x \in \text{supp}(\nu) $. For any integer $l \ge 0
$, we denote the open ball centered in $\mathbf x $ of radius
$\lambda^l $ by $ D_l = \{\mathbf y \in \Sigma: y_j = x_j \; \;
\forall \; 0 \le j < l \} $. Let $ H $ be a Walters module for the
potential $ A $. Given $ \epsilon > 0 $, we chose $ l $
sufficiently large (taking it larger than the number of
coordinates on which depends $\varphi $) in  such way that $
H(\lambda^l) < \epsilon/2 $.

As the constraint $ \varphi: \Sigma \to \mathbb Q^n $ is locally
constant, its image is reduced to a finite set of vectors with
rational coordinates. Suppose these numbers are written in
irreducible fractions and let $ Q
> 0 $ be the product of their denominators. In the same way, let
us consider $ q > 0 $ the product of the denominators of the
coordinates of $ r $.

The joint recurrence of $ \varphi $ assures there is a point $
\mathbf y \in \bigcap \Xi^j(D_l) \cap \text{\Large $\mathit b $}
(A) $. Then, we obtain a positive integer $ M_0 $ such that, for $
M \ge M_0 $, we have
$$\left | \frac{1}{M} S_M A(\mathbf y) - \int A \; d\nu \right | < \frac{\epsilon}{2}. $$
Besides, since in particular $\mathbf y \in \Xi^{M_0} (D_l) $, a
simple inductive argument gives positive integers $ L_1, \ldots,
L_{M_0} $ satisfying both $\sigma^{L_1 + \ldots + L_k} (\mathbf y)
\in \Xi^{M_0 - k} (D_l) $ and
$$\left \| S_{L_1 + \ldots + L_k} \varphi(\mathbf y) - (L_1 + \ldots + L_k) \varphi_*(\nu) \right \| < \frac{1}{qQ} \sum_{j = 1}^{k} {\frac{1}{2^{l+j}}} $$
for every $ k \in \{1, \ldots, M_0 \} $.

Put $ M = L_1 + \ldots + L_{M_0} \ge M_0 $. Take, then, the
periodic point $ \mathbf z \in \Sigma $ given by the repetition of
the word $ (y_0, \ldots, y_{M - 1}) $. Finally, let $\mu $ be the
$\sigma$-invariant probability measure by $\mathbf z $ defined. We
only need to verify that such periodic probability measure
accomplishes what is required.

Due to the fact we have taken $ l $ larger than the number of
coordinates on which depends the constraint $ \varphi $, we have $
\varphi(\sigma^j(\mathbf y)) = \varphi(\sigma^j(\mathbf z)) $ when
$ j \in \{0, \ldots, M - 1 \} $. Therefore,
$$ M \left \| \varphi_*(\mu) - r \right \| = \left \| S_M \varphi(\mathbf y) - M \varphi_*(\nu) \right \| < \frac{1}{qQ} \sum_{j = 1}^{M_0} {\frac{1}{2^{l+j}}} <
 \frac{1}{qQ} \cdot \frac{1}{2^l}. $$
Once $ QM\varphi_*(\mu) = Q S_M \varphi(\mathbf z) \in \mathbb Z^n
$, the inequality above assures $ \varphi_*(\mu) = r $. Besides,
observe that $ d_M (\mathbf y, \mathbf z) \le \lambda^l $ implies
$$\left | \int A \; d\mu - \frac{1}{M}S_M A(\mathbf y) \right | = \frac{1}{M} \left | S_M A(\mathbf z) - S_M A(\mathbf y) \right | \le \frac{1}{M} H(\lambda^l)
< \frac{\epsilon}{2}. $$ This ends the proof.
\end{proof}

There are two ways to interpret the conclusion of the theorem
above. The first is suggested by the well known fact according to
which a circle homeomorphism of rational rotation number possesses
a periodic point, whose period is equal to the denominator of the
rational number. Such point of view follows the same spirit, for
instance, of Franks theorem for certain rotation sets arising from
two-torus homeomorphisms homotopic to the identity (see
\cite{Franks}). In the context of the symbolic dynamics, a result
of this kind was obtained by Ziemian (see theorem 4.2 of
\cite{Ziemian}). The difference between the result of Ziemian and
the one obtained here is the transitivity hypothesis. We give up
this supposition, but we have introduced the condition of joint
recurrence. Thus, we have the immediate corollary.

\begin{5}
Suppose $ \varphi \in C^0(\Sigma, \mathbb Q^n) $ is a locally
constant function. Given $ r \in \varphi_*(\mathcal M_\sigma) \cap
\mathbb Q^n $, if there is in the fiber $ \varphi_*^{-1} (r) $ an
ergodic probability measure in relation to which $ \varphi $ is
joint recurrent, then in this fiber also exists a periodic
probability measure.
\end{5}

The second consequence of the theorem 9 is in the possibility of
supplying a special description to a beta function. We will
postpone the statement of the second corollary so that we can stop
shortly at this point.

In general, for an alpha function, we can indicate the
characterizations:
\begin{eqnarray*}
\alpha_{A, \varphi}(c) & = &
\min_{\mu \in \mathcal M_T} \int (\langle c, \varphi \rangle - A) \; d\mu \\
& = & \sup_{f \in C^0(X)} \min_{x \in X} (\langle c, \varphi \rangle - A + f - f \circ T) (x) \\
& = & \inf_{x \in \text{Reg}(\langle c, \varphi \rangle - A, T)}
\lim_{k \to \infty} \frac{1}{k} S_k(\langle c, \varphi \rangle - A)(x) \\
& = & \inf_{x \in X} \liminf_{k \to \infty}  \frac{1}{k}
S_k(\langle c, \varphi \rangle - A)(x),
\end{eqnarray*}
where $\text{Reg}(f, T) $ simply denotes the set of the points $ x
\in X $ for which is assured the existence of the limit of $
k^{-1} S_k f(x) $ when $ k $ tends to infinite. The first of the
equalities above, the reader will notice, comes directly from the
definition of the alpha function. The second expression is the
dual version of the previous one and it has been obtained recently
by Radu (see \cite{Radu}). Starting from the first, the last two
identities can be assured via Birkhoff's ergodic theorem. It is
possible to obtain these identities adapting lemmas contained in
the work of Hunt and Yuan (see the lemmas 2.3 and 2.4 of
\cite{HY}).

In relation to the representation of a beta function, we always
verify the dual formula
$$\beta_{A, \varphi}(h) = \inf_{(f,c) \in C^0(X) \times \mathbb R^n} \max_{x \in X} (A + f - f \circ T - \langle c, \varphi - h \rangle)(x). $$
When using the theorem of duality of Fenchel-Rockafellar, Radu
established such equality in \cite{Radu}. Starting from this
equality and using the other identities above, we get the
characterizations:
\begin{eqnarray*}
\beta_{A, \varphi}(h) & = &
\inf_{c \in \mathbb R^n} \beta_{A - \langle c, \varphi - h \rangle} \\
& = & - \sup_{c \in \mathbb R^n} \alpha_{A,\varphi - h}(c) \\
& = & \inf_{c \in \mathbb R^n} \sup_{x \in \text{Reg}(A - \langle
c, \varphi - h \rangle, T)}
\lim_{k \to \infty} \frac{1}{k} S_k(A - \langle c, \varphi - h \rangle)(x) \\
& = & \inf_{c \in \mathbb R^n} \sup_{x \in X} \limsup_{k \to
\infty}  \frac{1}{k} S_k(A - \langle c, \varphi - h \rangle)(x).
\end{eqnarray*}

Theorem 9 assures the following for subshifts of finite type.

\begin{5}
Let $ \varphi \in C^0(\Sigma, \mathbb Q^n) $ be a locally constant
constraint and $ A $ be a Walters potential. Taking $ r \in
\varphi_*(\mathcal M_\sigma) \cap \mathbb Q^n $, assume the
existence of an ergodic $(A,r)$-maximizing probability in relation
to which $ \varphi $ is joint recurrent. Then
$$ \beta_{A, \varphi}(r) = \sup \left \{\int A \; d\mu: \mu \in \varphi_*^{-1} (r), \mu \text{ periodic probability measure} \right \}. $$
\end{5}

The applicability of the corollary above is limited because the
need of finding a maximizing probability in relation to which the
constraint is joint recurrent. However, if we just concentrate on
constraints taking values in $ \mathbb Q $, the proposition 9, as
we saw, assures the recurrence.

The ergodic requirement also circumscribes the applicability of
corollary 12. Though, if we assume from now on in this section
that ours subshift of finite type $ \sigma: \Sigma \to \Sigma $ is
transitive, we have the following result.

\begin{4}
Let $ g \in C^0(\Sigma, \mathbb R^n) $ be a locally constant
function. Every point of the interior of the rotation set $
g_*(\mathcal M_\sigma) $ is a rotation vector of an ergodic
probability measure.
\end{4}

This theorem was obtained (see  theorem 4.6 of \cite{Ziemian}) for
Ziemian when the function $ g $ depends on two coordinates. By
passing to a higher block presentation of $ \Sigma $, the reader
familiarized with this argument will notice that the general case
is reduced to the situation treated by Ziemian (see also the
initial chapter of \cite{LM}).

We point out here the crucial importance of the locally constant
assumption. Let $ D \subset \mathbb R^n $ be a dense subset.
Consider $ \mathcal F_j (D) $ the collection of the functions of $
C^0 (\Sigma, D) $ that depend on $ j + 1 $ coordinates. It is
easily verified that $ \bigcup \mathcal F_j (D) $ is dense in $
C^0 (\Sigma, \mathbb R^n) $. Given a function $ g \in C^0 (\Sigma,
\mathbb R^n) $, we can without difficulty find a sequence $ \{g_j
\} $ convergent to $ g $ such that $ g_j \in \mathcal F_j (D) $
for each index $ j \ge 0 $.

\begin{2}
Suppose $\varphi \in C^0(\Sigma, \mathbb Q) $ is a locally
constant constraint and $ A $ is a Walters potential. Consider a
rational number $ r \in \text{int}(\varphi_*(\mathcal M_\sigma))
$. Then
$$ \beta_{A, \varphi}(r) = \sup \left \{\int A \; d\mu: \mu \in \varphi_*^{-1} (r), \mu \text{ periodic probability measure} \right \}. $$
\end{2}

\begin{proof}
Taking into account the theorem 10, fixed $ \epsilon > 0 $, it is
enough to assure the existence of an ergodic probability measure
$\nu $ with rotation number $ r $ satisfying $ \beta_{A,
\varphi}(r) - \epsilon < {\displaystyle \int A \; d\nu} $. Putting
$ \Phi = (\varphi, A) $, the strategy then lives in using the fact
of the graph of the application $ \beta_{A, \varphi} $ to do part
of the boundary of the rotation set $ \Phi_*(\mathcal M_{\sigma})
$. Thus, if this rotation set consists of a segment, the existence
of an ergodic probability measure as demanded follows  from
theorem 13.

It remains, therefore, to examine the other possibility: $
\text{int}(\Phi_*(\mathcal M_\sigma)) \ne \emptyset $. First, let
$ \{A_j \} \subset C^0(\Sigma) $ be a sequence convergent to $ A $
such that each function $ A_j $ depends on $ j + 1 $ coordinates.
Take any $ \eta > \beta_{A, \varphi}(r) - \epsilon/2 $ with $ (r,
\eta) \in \text{int}(\Phi_*(\mathcal M_\sigma)) $. When we put $
\Phi_j = (\varphi, A_j) $, from the proposition 4 it results $ (r,
\eta) \in \text{int}((\Phi_j)_*(\mathcal M_\sigma)) $ for an index
$ j $ sufficiently large, which can be supposed accomplishing
besides $ \| A_j - A \| _0 < \epsilon/2 $. By the theorem 13,
there exists an ergodic probability measure $ \nu \in \mathcal
M_\sigma $ satisfying $ (\Phi_j)_*(\nu) = (r, \eta) $, or better,
such that $ \varphi_*(\nu) = r $ and $ {\displaystyle \int A_j \;
d\nu} = \eta > \beta_{A, \varphi}(r) - \epsilon/2 $. However, once
$$ \left | \int A_j \; d\nu - \int A \; d\nu \right | \le \| A_j - A \| _0 < \frac{\epsilon}{2}, $$
it happens $ {\displaystyle \int A \; d\nu > \beta_{A, \varphi}(r)
- \epsilon} $.
\end{proof}

A natural question is: when $ \text{int}(\varphi_*(\mathcal
M_\sigma)) = \emptyset $?

In the setting we were in, there is a satisfactory answer. To
present it, though, it is convenient to detail a few more
properties of the general setting. A function $ g \in C^0(X) $ is
a (topological) cobord when there exists a function $ f \in C^0(X)
$ such that $ g = f \circ T - f $. Note that trivially every
cobord is a Walters function. Besides, two applications belonging
to $ C^0(X) $ are said cohomologous if their difference is a
cobord.

From results obtained by Bousch (in \cite{Bousch2}, consider
theorem 4 using theorem 1), it follows a particularly interesting
version of Liv\v{s}ic's theorem: an application $ f \in
C^0(\Sigma) $ is cohomologous to a constant if, and only if, $ f $
is a Walters function and $ \text{int}(f_*(\mathcal M_\sigma)) =
\emptyset $. A function locally constant, it is important to point
out, is a special example of Walters function.

\begin{5}
Let $ \varphi \in C^0(\Sigma, \mathbb Q) $ be a locally constant
constraint, not cohomologous to a constant. Assume $ A $ is a
Walters potential. For each $ c \in \mathbb R $, given $ \epsilon
> 0 $, there exist a rational number $ r \in
\text{int}(\varphi_*(\mathcal M_\sigma)) $ and a periodic
probability measure $ \mu \in \varphi_*^{-1} (r) $ satisfying $ c
r - {\displaystyle \int A \; d\mu} < \alpha_{A, \varphi}(c) +
\epsilon $.
\end{5}

\vspace{.8cm}

5. \textsc{{\large Sub-actions and Differentiability of Alpha
Functions}}

\vspace{.4cm}

We will obtain in the present section a result relating the
asymptotic behavior of optimal trajectories of certain sub-actions
and the differential of an alpha function. Let us recall that,
given a potential $ A $, an application $ u \in C^0(X) $ is a
sub-action (for $ A $) if
$$ A + u - u \circ T \le \beta_A. $$
General properties of sub-actions in distinct settings can be
found, for instance, in \cite{Branco, Bousch1, Bousch2, CLT,
Jenkinson3, LRR, LT1, LT2, PS, Savchenko, Souza}.

We denote
$$ \text{\Large $\mathit m$}_A =
\left \{\mu \in \mathcal M_T: \int A \; d\mu = \beta_A  \right \},
$$
the set of $A$-maximizing probabilities. Given a sub-action $ u $,
consider $ A^u = A + u - u \circ T $. It is easy to see that
$$ \text{\Large $\mathit m$}_A = \left \{\mu \in \mathcal M_T: \text{supp}(\mu) \subset (A^u)^{-1} \left ( \beta_A \right) \right \}. $$
Therefore, sub-actions help to locate the support of maximizing
probabilities. The compact set $ \mathbb M_A(u) =
(A^u)^{-1}(\beta_A) $ will be called the contact locus of the
sub-action $ u $\footnote{In \cite{CLT}, it was suggested to call
this set a Ma\~n\'e set.}. This is the set of points where the
above sub-action inequality turns out to be an equality.

We will consider $ (X, T) $ a transitive, expansive dynamical
system with a locally constant number of pre-images. Remind that
expansiveness means there exist $ \zeta > 0 $ and $ \kappa > 1 $
such that, if $ d(x, y) < \zeta $, then $ d(x, y) \kappa \le
d(T(x), T(y)) $. Besides, since the number of pre-images is
supposed locally constant, there is $ \xi > 0 $ such that,
whenever $ d(x', y') < \xi $ and $ x \in T^{-1}(x') $, we can find
$ y \in T^{-1}(y') $ accomplishing $ d(x, y) < \zeta $.

For a function $\theta$-H\"older $f $, the H\"older constant is
$$\text{H\"old}_\theta (f) = \sup_{d(x,y)>0} \frac{|f(x) - f(y) |}
{d(x,y)^\theta}. $$ As usual, we denote $ C^\theta (X) $ the
Banach space of $\theta$-H\"older functions with the norm $\|
\cdot \| _\theta = \text{H\"old}_\theta (\cdot) + \| \cdot \| _0
$.

We can now present a result that indicates how the variation of
the potential affects the sub-actions.

\begin{2}
Consider $ (X, T) $ a transitive, expansive dynamical system with
a locally constant number of pre-images. Let $\{B_j \} $ be a
sequence of $\theta$-H\"older functions converging in $C^\theta
(X) $ to a potential $ A $. Then, for each index $ j $, we can
find a sub-action $ v_j $ for the potential $ B_j $, so that any
accumulation point of the sequence $ \{v_j \} $ is a sub-action
for $ A $.
\end{2}

\begin{proof}

We have to show the existence of an equicontinuous and uniformly
bounded  sequence $\{v_j \} $. Indeed, as
$$\beta_{B_j} \ge B_j + v_j - v_j \circ T, $$
if $ u \in C^0(X) $ is an accumulation point of $ \{v_j \} $,
taking limit in $ j $, we immediately see that the function $ u $
is a sub-action for $ A $.

Given a  $\theta$-H\"older potential $ B $, it is possible to
obtain a sub-action $v $ for $B $ that satisfies
$$| v(x) - v(y) | \le \frac{\text{H\"old}_\theta(B)}
{\kappa^\theta - 1}\,\, d(x, y)^\theta,  \,\,\text {if}\,\,  d(x,
y) < \xi , \;\,\, \text {and also} $$
$$\| v \| _0 \le \text{H\"old}_\theta(B)
\left (\frac{2\xi^\theta}{\kappa^\theta - 1} + K
\text{diam}(X)^\theta \right), $$ being the positive integer $K
$\, depending just of $\xi $. For a proof of this statement, see
the reasoning of theorem 4.7 in \cite{Jenkinson4}.

As we are considering convergence in $ C^\theta (X) $, we clearly
obtain a sequence $ \{ v_j \} $ which is equicontinuous and
uniformly bounded.
\end{proof}

A sub-action $ u $ for a potential $\theta$-H\"older $ A $
satisfies
$$ u(x) - 2 u(T(x)) + u(T^2(x)) \ge - \text{H\"old}_\theta (A) \; d(x, T(x)) ^ \theta $$
for every $ x \in \mathbb M_A(u) $. Indeed, as
$$ (A + u - u \circ T)(x) = \beta_A $$
and
$$ (A + u - u \circ T)(T(x)) \le \beta_A, $$
we show the claim by simple subtraction. Besides, for a point $ x
$ belonging to the support of an $A$-maximizing probability, we
have
$$ \left | u(x) - 2 u(T(x)) + u(T^2(x)) \right | \le \text{H\"old}_\theta (A) \; d(x, T(x))^\theta. $$

Given a Walters potential $ A $, it is known the existence of a
sub-action $ u $ such that
$$ u(y) = \max_{T(x) = y} (A + u - \beta_A) (x). $$
This application $ u $ is called a calibrated sub-action for $ A
$.

We will suppose now a weaker assumption. We will consider a
transitive dynamical system $ (X, T) $ verifying the property of
weak expansion, that is, $ T^{-1}: \mathcal K(X) \to \mathcal K(X)
$ is 1-Lipschitz with respect to the Hausdorff metric. We can
assure the existence of calibrated sub-actions also in this
context (see \cite{Bousch2}).

For a calibrated sub-action $ u $, we will say that a sequence $
\{x_j \} \subset X $ is an optimal trajectory (associated to the
potential $ A $) when $ T(x_{j + 1}) = x_j $ and
$$ u(x_j) = A(x_{j + 1}) + u(x_{j + 1}) - \beta_A. $$

In general, as we remarked in last section, it is verified the
equality $ \alpha_{A, \varphi} (c) = - \beta_{A - \langle c,
\varphi \rangle} $. This is the last requirement for the
formulation of the next theorem.

\begin{4}
Let $ (X, T) $ be a transitive dynamical system satisfying the
property of weak expansion. Consider a  Walters potential $A \in
C^0(X) $, as well as  a Walters constraint $ \varphi \in C^0(X,
\mathbb R^n) $. Given an optimal trajectory $ \{x_j \} \subset X $
associated to the potential $A - \langle c, \varphi \rangle $, if
$\alpha_{A, \varphi} $ is differentiable at $c \in \mathbb R^n $,
we verify
$$\lim_{k \to \infty} \frac{1}{k} \sum_{j = 0}^{k - 1}
\varphi(x_j) = D\alpha_{A, \varphi}(c) .$$
\end{4}

\begin{proof}
Let $ u_c \in C^0(X) $ be the calibrated sub-action used in the
defitinion of the optimal trajectory $ \{x_j \} $. So we have
$$u_c (x_0) = u_c (x_k) + \sum_{j =
0}^{k - 1} \left [A(x_j) - \langle c, \varphi (x_j) \rangle +
\alpha_{A, \varphi}(c) \right]. $$ Consider $\rho > 0 $ and
$\gamma \in \mathbb R^n $ with $\| \gamma \| = 1 $. Taking any
calibrated sub-action $ u_{c + \rho \gamma} \in C^0(X) $ for the
potential $ A - \langle c + \rho \gamma, \varphi \rangle $, we
obtain
$$ u_{c + \rho \gamma} (x_0) \ge u_{c + \rho \gamma} (x_k) +
\sum_{j = 0}^{k - 1} \left [A(x_j) - \langle c + \rho \gamma,
\varphi (x_j) \rangle + \alpha_{A, \varphi}(c + \rho \gamma)
\right]. $$

From a simple subtraction, we get
\begin{eqnarray*}
- 2 \| u_c - u_{c + \rho \gamma} \|_0 & \le &
\sum_{j = 0}^{k - 1} \left[\langle\rho\gamma,\varphi(x_j)\rangle + \alpha_{A,\varphi}(c) - \alpha_{A,\varphi}(c+\rho\gamma)\right] \\
& = & \rho \left\langle \sum_{j = 0}^{k - 1} \varphi(x_j) - k
D\alpha_{A, \varphi}(c), \gamma \right\rangle + o(k\rho),
\end{eqnarray*}
therefore
$$\rho \left\langle \frac{1}{k}\sum_{j = 0}^{k - 1}
\varphi(x_j) - D\alpha_{A, \varphi}(c), \gamma \right\rangle =
O\left(\frac{1}{k}\right) + o(\rho). $$

Now taking limsup when $ k $ tends to infinite and using the fact
that $\rho $ can be arbitrarily small, we obtain $$\left\langle
\limsup_{k \to \infty}\frac{1}{k}\sum_{j = 0}^{k - 1} \varphi(x_j)
- D\alpha_{A, \varphi}(c), \gamma \right\rangle = 0 $$for all
$\gamma \in \mathbb R^n $ with $\| \gamma \| = 1 $, in other
words,
$$\limsup_{k \to \infty}\frac{1}{k}\sum_{j = 0}^{k - 1}
\varphi(x_j) = D\alpha_{A, \varphi}(c). $$

An analogous argument can be applied for the liminf and this
proves the theorem.
\end{proof}

\vspace {.4cm}

A similar result for the discrete Aubry-Mather problem is
presented in theorem 6.2 of \cite{Gomes}.


\begin{thebibliography}{99}

\bibitem{Atkinson}
G. Atkinson, Recurrence of co-cycles and random walks, \emph{The
Journal of the London Mathematical Society} \textbf{13} (1976),
486-488.

\bibitem{Branco}
F. M. Branco, \emph{Suba\c{c}\~ao para transforma\c{c}\~oes
unidimensionais}, PhD thesis, Universidade Federal do Rio Grande
do Sul, 2003.

\bibitem{Bousch1}
T. Bousch, Le poisson n'a pas d'ar\^etes, \emph{Annales de
l'Institut Henri Poincar\'e, Probabilit\'es et Statistiques}
\textbf{36} (2000), 489-508.

\bibitem{Bousch2}
T. Bousch, La condition de Walters, \emph{Annales Scientifiques de
l'\'Ecole Normale Sup\'erieure} \textbf{34} (2001), 287-311.

\bibitem{CI}
G. Contreras, R. Iturriaga, \emph{Global minimizers of autonomous
La-grangians}, 22$^\circ$ Col\'oquio Brasileiro de Matem\'atica,
IMPA, 1999.

\bibitem{CLT}
G. Contreras, A. O. Lopes, P. Thieullen, Lyapunov minimizing
measures for expanding maps of the circle, \emph{Ergodic Theory
and Dynamical Systems} \textbf{21} (2001), 1379-1409.

\bibitem{Franks}
J. Franks, Realizing rotation vectors for torus homeomorphisms,
\emph{Trans-actions of the American Mathematical Society}
\textbf{311} (1989), 107-115.

\bibitem{Garibaldi}
E. Garibaldi, \emph{Otimiza\c{c}\~ao erg\'odica: da
maximiza\c{c}\~ao relativa aos homeomorfismos expansivos}, PhD
thesis, Universidade Federal do Rio Grande do Sul, 2006.

\bibitem{Gomes}
D. A. Gomes, Viscosity solution method and the discrete
Aubry-Mather problem, \emph{Discrete and Continuous Dynamical
Systems, Series A} \textbf{13} (2005), 103-116.

\bibitem{Herman}
M. R. Herman, In\'egalit\'es a priori pour des tores lagrangiens
invariants par des diff\'eomorphismes symplectiques,
\emph{Publications Math\'ematiques de l'Institut des Hautes
\'Etudes Scientifiques} \textbf{70} (1989), 47-101.

\bibitem{HY}
B. R. Hunt, G. C. Yuan, Optimal orbits of hyperbolic systems,
\emph{Nonlin-earity} \textbf{12} (1999), 1207-1224.

\bibitem{Jenkinson1}
O. Jenkinson, \emph{Conjugacy rigidity, cohomological triviality
and barycentres of invariant measures}, PhD thesis, Warwick
University, 1996.

\bibitem{Jenkinson2}
O. Jenkinson, Geometric barycentres of invariant measures for
circle maps, \emph{Ergodic Theory and Dynamical Systems}
\textbf{21} (2001), 511-532.

\bibitem{Jenkinson3}
O. Jenkinson, Rotation, entropy, and equilibrium states,
\emph{Transactions of the American Mathematical Society}
\textbf{353} (2001), 3713-3739.

\bibitem{Jenkinson4}
O. Jenkinson, Ergodic optimization, \emph{Discrete and Continuous
Dynam-ical Systems, Series A} \textbf{15} (2006), 197-224.

\bibitem{Kwapisz1}
J. Kwapisz, Every convex polygon with rational vertices is a
rotation set, \emph{Ergodic Theory and Dynamical Systems}
\textbf{12} (1992), 333-339.

\bibitem{Kwapisz2}
J. Kwapisz, A toral diffeomorphism with a nonpolygonal rotation
set, \emph{Nonlinearity} \textbf{8} (1995), 461-476.

\bibitem{Kwapisz3}
J. Kwapisz, A priori degeneracy of one-dimensional rotation sets
for periodic point free torus maps, \emph{Transactions of the
American Mathematical Society} \textbf{354} (2002), 2865-2895.

\bibitem{LM}
D. Lind, B. Marcus, \emph{An introduction to symbolic dynamics and
coding}, Cambridge University Press, 1995.

\bibitem{LRR}
A. O. Lopes, V. Rosas, R. O. Ruggiero, Cohomology and
subcohomology for expansive geodesic flows, to appear in
\emph{Discrete and Continuous Dynamical Systems}.

\bibitem{LT1}
A. O. Lopes, P. Thieullen, Sub-actions for Anosov diffeomorfisms,
\emph{Ast\'erisque} \textbf{287} (2003), 135-146.

\bibitem{LT2}
A. O. Lopes, P. Thieullen, Sub-actions for Anosov flows,
\emph{Ergodic Theory and Dynamical Systems} \textbf{25} (2005),
605-628.

\bibitem{LT3}
A. O. Lopes, P. Thieullen, Mather measures and the Bowen-Series
transformation, to appear in \emph{Annales de l'Institut Henri
Poincare, Nonlinear Analysis}.

\bibitem{Mane}
R. Ma\~n\'e, Generic properties and problems of minimizing
measures of Lagrangian systems, \emph{Nonlinearity} \textbf{9}
(1996), 273-310.

\bibitem{MZ}
M. Misiurewicz, K. Ziemian, Rotation sets for maps of tori,
\emph{The Journal of the London Mathematical Society} \textbf{40}
(1989), 490-506.

\bibitem{PS}
M. Pollicott, R. Sharp, Livsic theorems, maximizing measures and
the stable norm, \emph{Dynamical Systems} \textbf{19} (2004),
75-88.

\bibitem{Radu}
L. Radu, Duality in thermodynamic formalism, \emph{preprint}.

\bibitem{Savchenko}
S. V. Savchenko, Cohomological inequalities for finite Markov
chains, \emph{Functional Analysis and Its Applications}
\textbf{33} (1999), 236-238.

\bibitem{Souza}
R. R. Souza, Sub-actions for weakly hyperbolic one-dimensional
systems, \emph{Dynamical Systems} \textbf{18} (2003), 165-179.

\bibitem{Ziemian}
K. Ziemian, Rotation sets for subshifts of finite type,
\emph{Fundamenta Mathematicae} \textbf{146} (1995), 189-201.

\end{thebibliography}
\end{document}